\documentclass[12pt,a4paper]{article}
\usepackage{epsfig} 
\usepackage{float}
\usepackage{subfigure}
\usepackage{graphicx}
\usepackage{amsmath}
\usepackage{amssymb}
\usepackage{setspace}
\usepackage{lineno}

\begin{document}
\title{On the Calculus of Functions with the Cantor-Tartan Support}
\author{   Alireza Khalili Golmankhaneh $^{a*}$}
\date{}
\maketitle \vspace{-9mm}
\begin{center}
$^a$ Department of Physics, Urmia Branch, Islamic Azad University, Urmia, Iran\\
 \emph{*E-mail address}: a.khalili@iaurmia.ac.ir
 \end{center}

\begin{abstract}
\doublespacing
In this manuscript, integral and derivative of the functions with  Cantor-Tartan spaces is defined. The generalization of  standard calculus which is called $F^{\alpha}$-calculus utilized to obtain the integral and derivative of the functions on the Cantor-Tartan with different dimensions. Differential equation
involving the new derivatives are solved.  The illustrative examples are used to present the details.
\vspace{.5cm} \noindent {\it {Keywords:}} $F^{\alpha}$-calculus;  Staircase function ; Cantor-Tartan support; Fractional differential equation\\

\end{abstract}
The fractal shapes and objects are seen in the nature, e.g.,  clouds,  mountains,  coastlines, human body, and etc. The geometry of the fractals were studied   \cite{book-1}. The analysis on fractals were established, using different methods, such as fractional calculus, probability theory, measure theory,  fractional spaces, and time scale theory  by many researchers and found many applications \cite{book-2,book-3,book-3,book-4,book-5,book-6,book-7,book-8,book-88,book-9,book-11}.~The fractional derivatives  have non-local property which are suitable to model the process with the memory effect, non-conservative systems \cite{book-12,book-13,book-14,book-15,book-16,book-17,book-18,book-19,book-20,book-21}.~The fractional calculus, which involves derivatives with arbitrary orders, has applied on the process with  fractals structures  \cite{conection1,conection2,conection3}. The anomalous diffusion on fractals was formulated which included sub- and supper diffusion in view of different random walks \cite{d-5,d-6,d-7}.~Fractal antennas are small but have wide-band radiations which make them useful in microwave communications \cite{book-444,book-4445}. Laminar flow of a fractal fluid in a cylindrical tube was studied using the homogeneous flow in a fractional dimensional space \cite{Paper-new2mn3}.~On the Cantor cubes, the Maxwell's equation  were obtained and, as an application, the electric field due to Cantor dust was obtained \cite{Paper-new2mn}. Recently, $F^{\alpha}$-calculus which is algorithmic was suggest and applied for modeling some  physical process \cite{Gangal-1,Gangal-2}.  As a pursuit of  these researches,  $F^{\alpha}$-calculus is generalized and utilized in optics  and mechanics \cite{Golmankhaneh-1,Golmankhaneh-2,Golmankhaneh-3,Ali-1}. The non-local integrals and derivatives are defined  on the Cantor set which are used to model the fractal ideal solids and the fractal ideal fluids \cite{Golmankhaneh-4,Golmankhaneh-5}.\\
This manuscript is built as follows:\\
In Section  \ref{2-dim} we set up the notation and terminology of $F^{\alpha}$-calculus on the fractals that imbedded  in $\Re^{2}$ \cite{Paper-new2c3,Paper-new2c4}. In Section  \ref{3-dim} we  study  some examples on the Cantor-Tartan using suggested definitions. Section  \ref{4-dim} is devoted to the conclusion.

\section{Terminology and notations \label{2-dim}}
\begin{figure}[H]
     \includegraphics[scale=0.5]{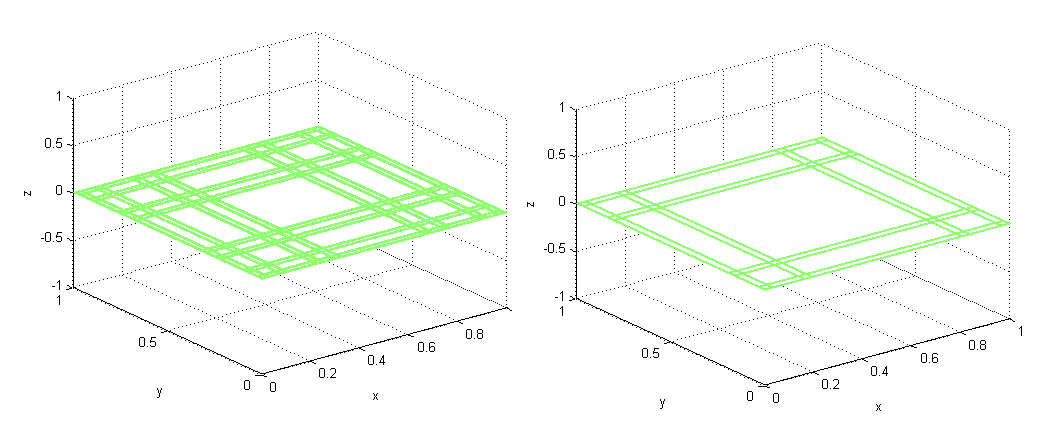}
    \caption{ The Cantor-Tartan with different dimensions are plotted.}
    \label{fig:1}
\end{figure}

In the section, we review and define the basic tools we need in our work.\\
Fractals are the sets with the self-similar properties such that  theirs fractal dimension exceeds from theirs topological dimension. The calculus on the Cantor-Tartan $\mathfrak{F}=F\cup F\subset \Re^2$  where $F=C\times \Re$, $C$ is Cantor set and $\Re$ is real line \cite{book-2}.~Let  $\mathfrak{F}$ be Cantor-Tartan that is subset of $I=[a,b] \times [c,d],~a,~b,~c,~d \in\Re$ (Real-line).~We sketch in Figures [\ref{fig:1}]  the  Cantor-Tartan with different fractal box or similar  dimensions.~The Cantor-Tartan is established by Cartesian product of Cantor set and real line (F) and union of them. Here, we consider  finite iteration which is approximation of fractals spaces.

~The flag function for the Cantor-Tartan $\mathfrak{F}=F\cup F\subset \Re^2$, which is denoted by   $\Theta(\mathfrak{F},I)$,  is defined  as follows
\begin{equation}
\Theta(\mathfrak{F},I)=\begin{cases}
    1, ~~~\textmd{if}~~~ \mathfrak{F}\cap I \neq \emptyset,\\
    0, ~~~~\textmd{otherwise}.
\end{cases}
\end{equation}
The subdivision $P_{[a,b]}$ is
\begin{equation}
P_{[a,b]\times [c,d]}=\{x_{0}=a, x_{1}, x_{2},...,x_{n}=b\}\times \{y_{0}=c, y_{1}, y_{2}, ...,~y_{m}=d\},
\end{equation}
where $\times$ denotes the Cartesian product.\\
The $\gamma^{\zeta}(\mathfrak{F},a,b,c,d)$ is defined
\begin{align}
&\gamma^{\zeta}(\mathfrak{F},a,b,c,d)=\nonumber \\&
\lim_{\delta\rightarrow 0}
\inf_{P_{[a,b]\times [c,d]}:|P|\leq\delta} \sum_{j=1}^{n}\sum_{i=1}^{n}\frac{(x_{i}-x_{i-1})^{\alpha}}{\Gamma
(\alpha+1)}\frac{(y_{j}-y_{j-1})^{\beta}}{\Gamma
(\beta+1)}~\Theta(\mathfrak{F},[x_{i-1},x_{i}])~\Theta(\mathfrak{F},[y_{j-1},y_{j}]),
\end{align}
where $\zeta=\max\{\alpha+\beta\}$, $\alpha=\beta=1+Dim(C)$  and $|P|$ is

\begin{equation}\label{c9}
  |P|=\max_{1\leq i \leq n,~ 1\leq j \leq n} (x_{i}-x_{i-1})\times(x_{j}-x_{j-1}).
\end{equation}
The integral staircase function  $S_{\mathfrak{F}}^{\zeta}(x,y)$  for Cantor-Tartan $\mathfrak{F}$ is defined
\begin{equation}\label{t}
    S_{F}^{\zeta}(x,y)=\begin{cases}
    \gamma^{\zeta}(F,a_{0},c_{0},x,y), ~~~~\text{if}, ~~~~~x\geq a_{0},~~y\geq c_{0};\\
    -\gamma^{\zeta}(F,a_{0},c_{0},x,y), ~~~~\text{otherwise},
\end{cases}
\end{equation}
where $a_{0},~c_{0}$ are  arbitrary real numbers.

\begin{figure}[H]
   \includegraphics[scale=0.5]{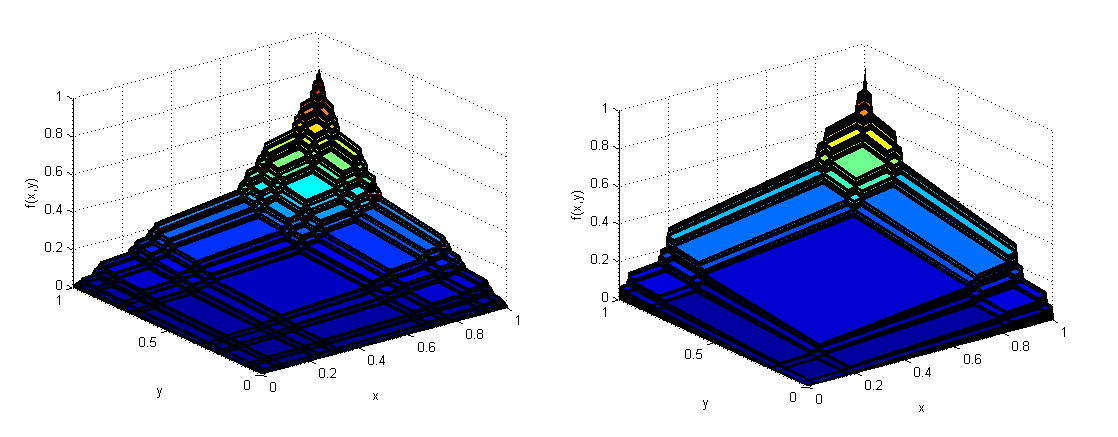}
    \caption{Staircase function corresponding to Cantor-Tartan with different dimensions are presented.}
    \label{fig:2}
\end{figure}

The integral staircase function for the different Cantor-Tartan with the different dimension  are plotted in  Figure [\ref{fig:2}].\\
The $\zeta$-dimension of $ \mathfrak{F}\cap [a,b]\times [c,d]$ is given by
\begin{eqnarray}\label{sa}
    \dim_{\zeta}(\mathfrak{F}\cap[a,b]\times[c,d])&=&\inf\{\zeta:\gamma^{\zeta}(\mathfrak{F},a,c,b,d) =0\},\\
    &=&\sup\{\zeta:\gamma^{\zeta}(\mathfrak{F},a,c,b,d)=\infty\}.\nonumber
\end{eqnarray}
A point $(x,y)$ is called a point of change of a function $f(x,y)$ if it is not constant over any open set $[a,b]\times [c,d]$ involving $(x,y)$. The set of all  points of change of $f(x,y)$ are indicated by $Sch f$.\\ The $Sch(S_{\mathfrak{F}}^{\zeta}(x,y))$ is called $\zeta$-perfect if $S_{\mathfrak{F}}^{\zeta}(x,y)$ is finite for all $(x,y)\in \Re$.\\
Let $f(x,y)$ be a  bounded function on $\mathfrak{F}$  then we define
\begin{eqnarray}
  M[f,\mathfrak{F},I] &=& \sup_{(x,y)\in \mathfrak{F}\cap I} f(x,y),~~~~\text{if}~~~\mathfrak{F}\cap I\neq 0. \\
  &=& 0, ~~~~~ \text{otherwise},\nonumber
\end{eqnarray}
and similarly
\begin{eqnarray}
  m[f,\mathfrak{F},I] &=& \inf_{(x,y)\in \mathfrak{F}\cap I} f(x,y),~~~~\text{if}~~~\mathfrak{F}\cap I\neq 0.\\
  &=& 0, ~~~~~ \text{otherwise}.\nonumber
\end{eqnarray}
Now, upper  $U^{\zeta}$-sum and  lower $L^{\zeta}$-sum for the
function $f(x,y)$ over the subdivision $P$ are given as follows
\begin{equation}\label{y}
    U^{\zeta}[f,\mathfrak{F},P]=\sum_{i=1}^{n}M[f,\mathfrak{F},
    [(x_{i-1},y_{i-1}),(x_{i},y_{i})]](S_{\mathfrak{F}}^{\zeta}(x_{i},y_{i})
-S_{\mathfrak{F}}^{\zeta}(x_{i-1},y_{i-1})),
\end{equation}
and
\begin{equation}\label{y}
    L^{\zeta}[f,\mathfrak{F},P]=\sum_{i=1}^{n}m[f,\mathfrak{F},
    [(x_{i-1},y_{i-1}),(x_{i},y_{i})]]
(S_{\mathfrak{F}}^{\zeta}(x_{i},y_{i})-S_{\mathfrak{F}}^{\zeta}(x_{i-1},y_{i-1})).
\end{equation}
The $f(x,y)$ is called $F^{\zeta}$-integrable on Cantor-Tartan $\mathfrak{F}$ if we have
\begin{equation}\label{cft}
    \underline{\int_{(a,c)}^{(b,d)}}f(x,y) d_{F}^{\alpha}x d_{F}^{\beta}y=\sup_{P_{[a,b]}} L^{\zeta}[f,\mathfrak{F},P]=
    \overline{\int_{(a,c)}^{(b,)d}}f(x,y) d_{F}^{\alpha}x d_{F}^{\beta}y=\inf_{P_{[a,b]}} L^{\zeta}[f,\mathfrak{F},P].
\end{equation}
The $F^{\zeta}$-integral is denoted by $\int_{(a,b)}^{(b,d}f(x,y) d_{F}^{\alpha}x d_{F}^{\beta}y$.\\
Let $\mathfrak{F}$ be $\zeta$-perfect set then  we define  $F^{\alpha}$-partial derivative  of $f(x,y)$ respect to $x$ as
\begin{equation}
    ^{x}D_{\mathfrak{F}}^{\alpha}f(x,y)=\begin{cases}
    \mathfrak{F}-\lim_{(x',y)\rightarrow (x,y)} \frac{f(x',y)-f(x,y)}{S_{\mathfrak{F}}^{\zeta}(x',y)-S_{\mathfrak{F}}^{\zeta}(x,y)}, ~~~\text{if} ~~~(x,y)\in \mathfrak{F},\\
    0, ~~~~~~~~~~~~~\textmd{otherwise}.
    \end{cases}
 \end{equation}
if the limit exists.
\section{The functions with Cantor-Tartan support \label{3-dim}}
The task is now to apply definitions on examples. In this section, we give some examples to show more details.\\
\subsection{Example 1.} 
\begin{figure}[H]
  \includegraphics[scale=0.5]{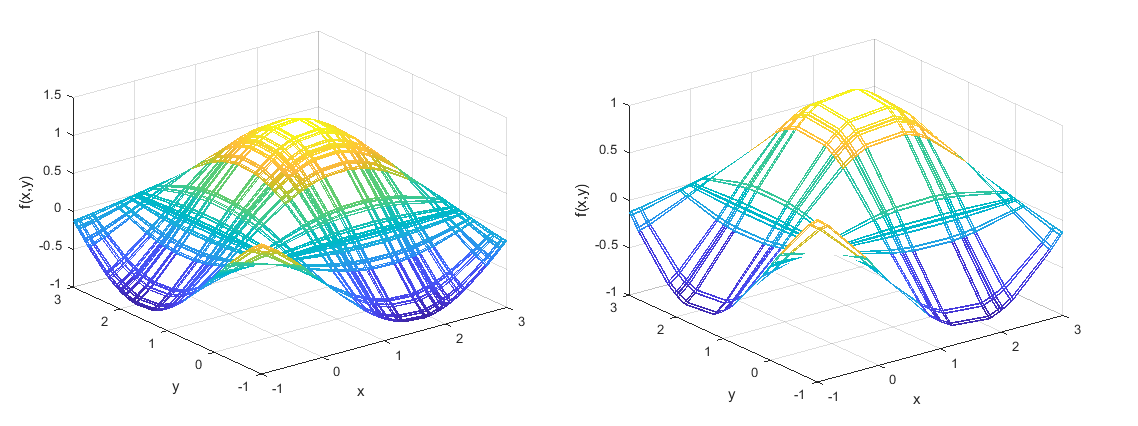}
    \caption{We plot the graph of function $\sin(\chi_{F}^{\alpha}x)\sin(\chi_{F}^{\beta}y)$ with the different Cantor-Tartan support.}
    \label{fig:3}
\end{figure}

Let us consider a function with Cantor-Tartan support  with the different fractal dimensions
\begin{equation}\label{f}
 f(x,y)=\sin(\chi_{F}^{\alpha}x)\sin(\chi_{F}^{\beta}y),~ (x,y)\in \mathfrak{F},
\end{equation}
 where $\chi_{F}^{\alpha},~\chi_{F}^{\beta}$ are characteristic function for fractal sets $F$ \cite{Golmankhaneh-5,Gangal-1}.~The graph of the $f(x,y)$ is shown in Figures [\ref{fig:3}].

  The fractal integral of $f(x,y)$ on Cantor-Tartan $\mathfrak{F}\subset[0,1]\times [0,1]$ is as follows:
\begin{align}
 &g(x,y)\bigg|_{(x=y=1)}=\int_{0}^{x}\int_{0}^{y}\sin(\chi_{F}^{\alpha}x)
 \sin(\chi_{F}^{\beta}y)d_{F}^{\alpha}x' d_{F}^{\beta}y' \bigg|_{(x=y=1)}= \nonumber \\ & \int_{0}^{1}\cos(S_{F}^{\beta}(y'))\sin(\chi_{F}^{\alpha}x)\bigg|_{0}^{1} d_{F}^{\alpha}x' = \nonumber \\ & \int_{0}^{1}[\cos(S_{F}^{\beta}(1))\sin(\chi_{F}^{\alpha}x)-\cos(S_{F}^{\beta}(0))\sin(\chi_{F}^{\alpha}x)]~d_{F}^{\alpha}x'= \nonumber\\& \cos(\Gamma(1+\beta))\cos(S_{F}^{\alpha}(x'))-\cos(S_{F}^{\alpha}(x'))\bigg|_{S_{F}^{\alpha}(0)}^{S_{F}^{\alpha}(1)}= \nonumber\\&
\cos(\Gamma(1+\beta))\cos(\Gamma(1+\alpha))-\cos(\Gamma(1+\alpha))-\cos(\Gamma(1+\beta))+1\nonumber\\&
 =\left\{
    \begin{array}{ll}
      0.13, & \alpha=\beta=0.6~~\textmd{With~support~Cantor-Tartan~dimesion~ 1.2} \\
      0.16, & \alpha=\beta=0.8~~\textmd{With~support~Cantor-Tartan dimesion ~1.6}
    \end{array}
  \right.
\end{align}
here $S_{F}^{\alpha}(1)=\Gamma(1+\alpha),~S_{F}^{\beta}(1)=\Gamma(1+\beta),~S_{F}^{\alpha}(0)=S_{F}^{\beta}(0)=0$ .
\begin{figure}[H]
\includegraphics[scale=0.5]{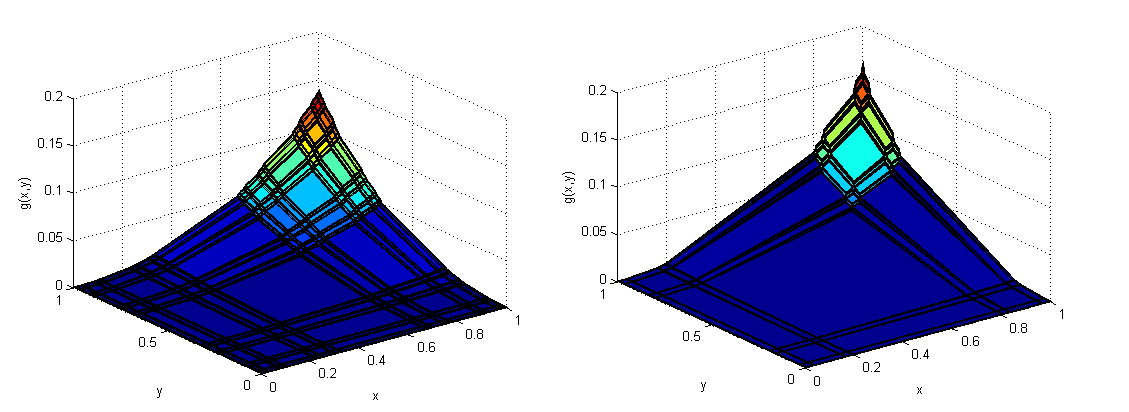}
    \caption{Fractal integral  of $\sin(\chi_{F}^{\alpha}x)\sin(\chi_{F}^{\beta}y)$ on the different  Cantor-Tartan are sketched . }
    \label{fig:4}
\end{figure}
In Figure [\ref{fig:4}] we sketch the $g(x,y)$ which is called fractal integral of $f(x,y)$.\\
\subsection{Example 2.} Consider a function with Cantor-Tartan support as follows
\begin{equation}\label{mn23}
f(x,y)=\sin(\chi_{F}^{\alpha}x+\chi_{F}^{\beta}y),~ (x,y)\in \mathfrak{F}
\end{equation}

\begin{figure}[H]
   \includegraphics[scale=0.5]{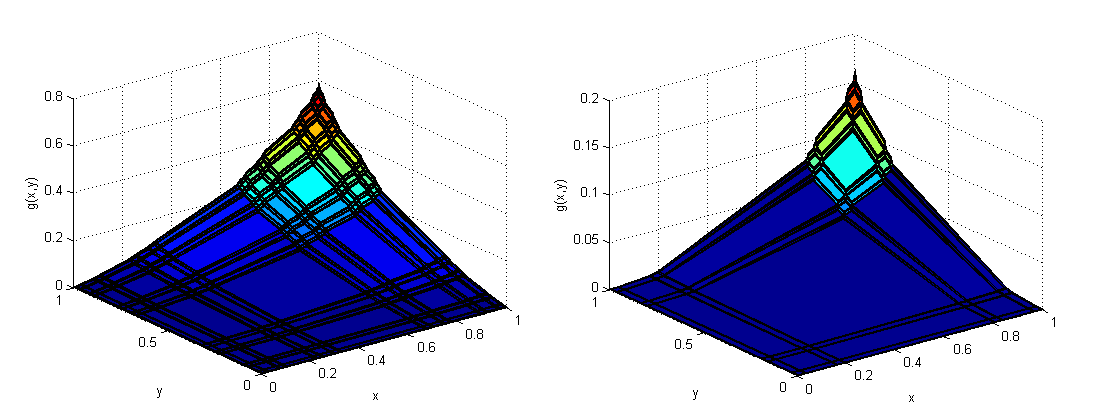}
    \caption{We plot the graph of function $\sin(\chi_{F}^{\alpha}x+\chi_{F}^{\beta}y)$ with the different Cantor-Tartan support. }
    \label{fig:5}
\end{figure}

We present in Figure [\ref{fig:5}] graph of Eq. (\ref{mn23}). The fractal integral of Eq. \eqref{mn23} is as follows
\begin{align}
 &g(x,y)\bigg|_{(x=y=1)}=\int_{0}^{x}\int_{0}^{y}\sin(\chi_{F}^{\alpha}x+\chi_{F}^{\beta}y)d_{F}^{\alpha}x'~ d_{F}^{\beta}y' \bigg|_{(x=y=1)}= \nonumber \\ & \int_{0}^{1}-\cos(\chi_{F}^{\alpha}x+S_{F}^{\beta}(y))\bigg|_{0}^{1} d_{F}^{\alpha}x' = \nonumber \\ & \int_{0}^{1}[-\cos(\chi_{F}^{\alpha}x+S_{F}^{\beta}(1))+\cos(\chi_{F}^{\alpha}x)]~d_{F}^{\alpha}x'= \nonumber\\& -\sin(S_{F}^{\alpha}(x')+\Gamma(1+\beta))+\sin(S_{F}^{\alpha}(x'))\bigg|_{0}^{1}= \nonumber\\&
 =-\sin(\Gamma(1+\alpha)+\Gamma(1+\beta))+\sin(\Gamma(1+\alpha))
+\sin(\Gamma(1+\beta))\nonumber\\&
 =\left\{
    \begin{array}{ll}
      0.581, & \alpha=\beta=0.6~~\textmd{With~support~Cantor-Tartan~dimesion~ 1.2} \\
      0.64, & \alpha=\beta=0.8~~\textmd{With~support~Cantor-Tartan dimesion ~1.6}
    \end{array}
  \right.
\end{align}
where $S_{F}^{\alpha}(1)=\Gamma(1+\alpha),~S_{F}^{\beta}(1)=\Gamma(1+\beta),~S_{F}^{\alpha}(0)=S_{F}^{\beta}(0)=0$.
\begin{figure}[H]
\includegraphics[scale=0.5]{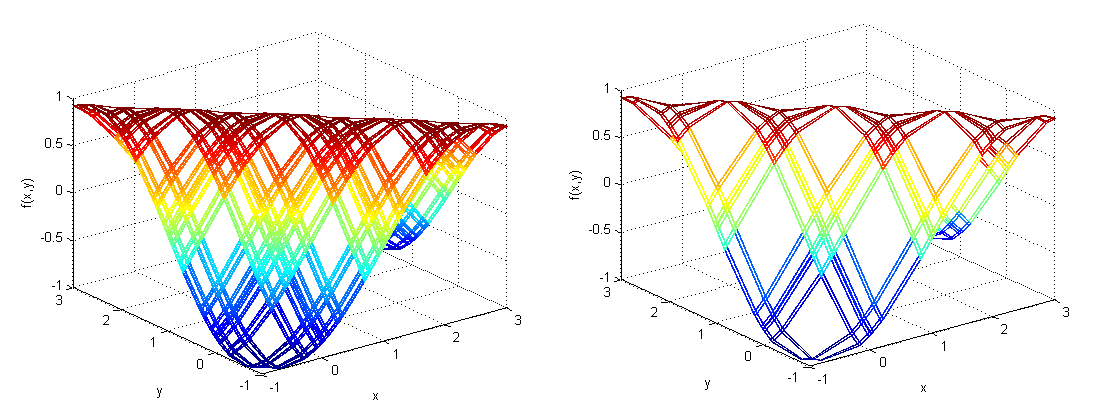}
\caption{Fractal integral  of $\sin(\chi_{F}^{\alpha}x+\chi_{F}^{\beta}y)$ on different Cantor-Tartan  are shown.}
    \label{fig:6}
\end{figure}
In Figure [\ref{fig:6}] we plot fractal integral of $\sin(\chi_{F}^{\alpha}x+\chi_{F}^{\beta}y)$ versus to different  Cantor-Tartan support with different dimensions.\\

\subsection{Example 3.}
Consider a function with  the Cantor-Tartan as follows
\begin{equation}\label{c}
f(x,y)=S_{F}^{\alpha}(x)^2+S_{F}^{\beta}(y)^2.
\end{equation}
The fractal partial derivatives of $f(x,y)$ respect to $x$ and $y$  are
\begin{equation}\label{8b}
  ^{x}D^{\alpha}_{F}f(x,y)=2 S_{F}^{\alpha}(x)\chi_{F},~~~^{y}D^{\alpha}_{F}f(x,y)=2 S_{F}^{\alpha}(y) \chi_{F}.
\end{equation}

\begin{figure}[H]
\includegraphics[scale=0.5]{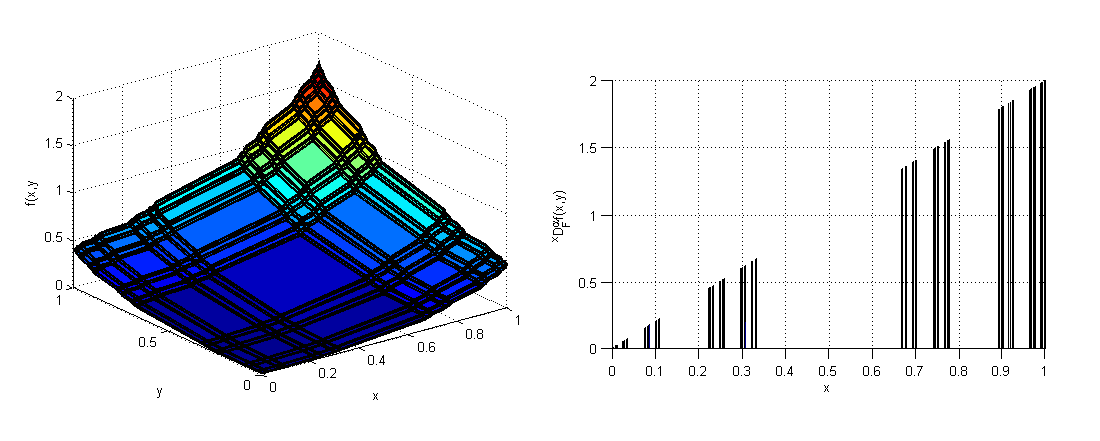}
      \caption{We have plotted $S_{F}^{\alpha}(x)^2+S_{F}^{\beta}(y)^2$ and  $^{x}D^{\alpha}_{F}(x,y)$.}
    \label{fig:7}
\end{figure}

In Figure [\ref{fig:7}] we plot Eq. (\ref{c}) and it's partial derivative which is shown $^{x}D^{\alpha}_{F}(x,y)$.

\textbf{Remark:} One can obtain standard results by choosing $\alpha=1,~\beta=1,~S=\Re$ which leads to $ S_{F}^{\alpha}(x)=x,~ S_{F}^{\beta}(y)=y$.

\section{Conclusion \label{4-dim}}
In this work, we define the local derivative and integral on Cantor-Tartan. The standard calculus can not be applied to integrate  and differentiate for the function on this fractals. Therefore, we need a new type of calculus to calculate the physical properties and describe  phenomena on fractals. As a result, the $F^{\alpha}$-calculus on Cantor-Tartan which has fractal dimension $1<\zeta<2$ is given.  More, we recall  the standard calculus results which shows that suggested definitions are the generalization of standard calculus. Three illustrative examples were investigated and the corresponding graphs of the functions were drawn.


\end{document}